\documentclass[12pt]{article}
\usepackage{amssymb,amsmath,amsthm}

\title{Painlev\'e Test and the Resolution of Singularities for Integrable Equations}
\author{Jishan Hu, 
Min Yan \\
The Hong Kong University of Science and Technology}

\newtheorem{theorem}{Theorem}

\newtheorem{proposition}[theorem]{Proposition}
\newtheorem*{theoremA*}{Theorem A}
\newtheorem*{theoremB*}{Theorem B}
\newtheorem*{corollary*}{Corollary}

\theoremstyle{definition}
\newtheorem{definition}[theorem]{Definition}
\newtheorem*{definition*}{Definition}
\newtheorem*{case*}{Case}
\newtheorem*{subcase*}{Subcase}

\theoremstyle{remark}

\numberwithin{equation}{section}

\begin{document}
\maketitle

\tableofcontents

\section{Introduction}

Since Kowalevskaya's monumental work \cite{kowa}, the Painlev\'e test has been the most widely used and the most successful technique for detecting integrable differential equations \cite{ars,cm,fikn}. The test has been applied to many differential equations and, for those passing the Painlev\'e test, the indicators of integrability such as the existence of sufficiently many conservation laws, the Lax pair, the Darboux transform, the B\"acklund transform have always been found. 

The Painlev\'e test itself is a recipe for finding all formal Laurent series solutions with movable singularities. An $n$-th order system {\em passes the Painlev\'e test} if such formal Laurent series solutions admit $n$ free parameters (including the location $t_0$ of the the movable singularity as one free parameter). Since the solutions of an $n$-th order system should have $n$ degrees of freedom, the Laurent series solutions, with maximal number of free parameters, should include all the solutions. Then one hopes that all movable singular solutions are poles and concludes that all solutions are single valued. This is the heuristic reason why passing the Painlev\'e test means integrability. Indeed such reason underlies the classification of integrable equations by Painlev\'e and the others \cite{pain}. 

In \cite{am1,am2}, Adler and van Moerbeke put the heuristic reason on solid foundation for the very nice case of algebraically completely integrable systems. They used the toric geometry to construct a complete phase space, and thereby gave a satisfactory explanation for the relation between the Painlev\'e test and the integrability. In \cite{es}, Ercolani and Siggia pointed out that algebraic geometry is not needed for integrable Hamiltonian systems. They suggested using the expansion of the Hamilton-Jacobi equation to construct the change of variable used for completing the phase space. They showed that the method often works through many examples, but did not prove that the method always works.

In this paper, we show that, not only is algebraic geometry not needed, the Hamiltonian set up by Ercolani and Siggia is also not needed. In fact, with virtually no condition, the existence of the change of variable needed for completing the phase space is equivalent to passing the Painlev\'e test.

\begin{theoremA*}
A regular system of ordinary differential equations passes the Painlev\'e test if and only if there is a triangular change of variable, such that the system is converted to another regular system, and the Laurent series solutions produced by the Painlev\'e test are converted to power series solutions.
\end{theoremA*}

For the exact meaning of ``passing the Painlev\'e test'', see the rather straightforward definition in the beginning of Section \ref{general}. Also see \cite{hy2} for some concrete examples. The condition is as straightforward and elementary as can be. Although the proof is also rather elementary, we believe there is some interesting underlying algebraic structure that is worth further exploration.  

The theorem gives a key connection between the two major works by Kowalevskaya \cite{hy5}. As in the theorem of Cauchy and Kowalevskaya, we consider a regular system of differential equations with complex analytic functions on the right side. Given a movable singular solution, represented by Laurent series in $(t-t_0)$ for variables $u_i$, changing $u_i$ to $u_i^{-1}$ certainly regularizes the variables. However, such simple regularization will create singular differential equations for the new variables. By the resolution of singularity, we mean a change of variable that regularizes both the solution and the equation. This is not always possible, and the condition for the resolvability of the singularity is exactly passing the Painlev\'e test.

The regularization is achieved by a {\em triangular change of variable}
\begin{align*}
u_1 
&= \tau^{-k}, \\
u_i 
&= a_i(t,\tau,\rho_2,\dots,\rho_{i-1})+b_i(t,\tau,\rho_2,\dots,\rho_{i-1})\rho_i,\quad 1<i\le n,
\end{align*}
where $-k$ is the leading order of the Laurent series solution for $u_1$ obtained in the Painlev\'e test, $a_i,b_i$ are meromorphic in $\tau$ and analytic in the other variables, and $b_i$ does not take zero value. Such change of variable can be easily inverted, at the cost of introducing one $k$-th root. Then it is easy to see that the Laurent series for $u_1,\dots,u_n$ correspond to a power series for $\tau$ satisfying $\tau(t_0)=0$, $\tau'(t_0)\ne 0$, and the Laurent series for $\rho_2,\dots,\rho_n$. In case the system for $u_1,\dots,u_n$ passes the Painlev\'e test, we can find suitable $a_i,b_i$ so that the Laurent series for $\rho_2,\dots,\rho_n$ are always power series. In fact, in our change of variable formula \eqref{change1} and \eqref{triangle} (or \eqref{change1} and \eqref{change3}), we get very specific forms for $a_i$ and $b_i$. Then applying Cauchy theorem to the regular system for $\tau,\rho_2,\dots,\rho_n$ and the convert to the original system, we conclude the following. See \cite[Section 3.8.6]{go} for more discussion on the convergence of formal Laurent series solutions.

\begin{corollary*}
If a regular system of ordinary differential equations passes the Painlev\'e test, then the Laurent series solutions obtained from the test always converge. 
\end{corollary*}

Section \ref{general} is devoted to the rigorous proof of Theorem A. We will be rather conservative in assuming that the right side of the differential equations are polynomials. The assumption is used to make sure that, after substituting in the Laurent series, the right side again becomes a Laurent series. Certainly we can make sense of this with functions other than polynomials, but the spirit of Theorem A is already fully reflected in the polynomial case. 

Ercolani and Siggia's work suggests that there should be a version of Theorem A that is compatible with Hamiltonian structure. 

\begin{theoremB*}
If a regular Hamiltonian system of ordinary differential equations passes the Painlev\'e test in the Hamiltonian way, then there is a canonical triangular change of variable, such that the system is converted to another regular Hamiltonian system, and the Laurent series solutions are converted to power series solutions. Moreover, if the system is autonomous, then the new Hamiltonian function is obtained by substituting the new variables. If the system is not autonomous, then the new Hamiltonian function is obtained by substituting the new variables and then dropping the singular terms.
\end{theoremB*}

The exact meaning of the ``Hamiltonian way'' is given in Section \ref{hamilton}. See \cite{hy3} for some concrete examples.

In view of the rather ad hoc nature of the Painlev\'e test, we further provide a rigorous foundation of the Painlev\'e test in Section \ref{foundation}. Our setup is by no means the broadest possible (and is therefore only ``a'' foundation), but is satisfied by most examples and is the most common way the conditions of our theorems are met. 

We emphasize that the results of this paper are local, in fact as local as the theorem of Cauchy and Kowalevskaya. Although the Painlev\'e test is local, it is supposed to test the global properties of solutions. The regularization in this paper may be used to further prove the global property for some specific integrable equations, but here we are not making any general statement about the global property. See \cite{hl1,hl2,hl3,hy6,st} for specific examples of global Painlev\'e analysis.

Finally, we remark that our theorems do not deal with the lower balances between the principal balances, as discussed in \cite{am1,am2,es} as part of the process of completing the phase space. If there are no lower balances, then our theorems effectively show that passing the Painlev\'e test implies the completion of the phase space and thereby explains the integrability. Further extensions of our theorems to lower balances are needed for the more general case.

\section{Resolution of Movable Singularity} 
\label{general}

Consider a system of ordinary differential equations
\[
u_i'=f_i(t,u_1,\dots,u_n),\quad 1\le i\le n.
\]
We restrict the discussion to $f_i(t,u_1,\dots,u_n)$ being analytic in $t$ and polynomial in $u_1,\dots,u_n$. At the end of this section, we will point out what is exactly needed for more general $f_i$.

A {\em balance} for the system is a formal Laurent series solution, with {\em movable singularity} $t_0$. The balance often admits several free parameters, similar to the initial values in formal power series solutions. Therefore a balance is typically given by ($r_1$ is reserved for the first free parameter $t_0$)
\begin{align}
u_i
&= a_{i,0}(t_0)(t-t_0)^{-k_i}+\dots+a_{i,\lambda_1-1}(t_0)(t-t_0)^{\lambda_1-1-k_i} \nonumber \\
&\quad +a_{i,\lambda_1}(t_0,r_2,\dots,r_{n_1})(t-t_0)^{\lambda_1-k_i} \nonumber \\
&\qquad +\dots+a_{i,\lambda_2-1}(t_0,r_2,\dots,r_{n_1})(t-t_0)^{\lambda_2-1-k_i} \nonumber \\
&\quad + \cdots  \label{balance0} \\
&\quad +a_{i,\lambda_s}(t_0,r_2,\dots,r_{n_s})(t-t_0)^{\lambda_s-k_i} \nonumber \\
&\qquad +\dots+a_{i,j}(t_0,r_2,\dots,r_{n_s})(t-t_0)^{j-k_i}+\cdots, \nonumber 
\end{align}
where $a_{i,j}$ is analytic in its variables. The free parameters are called {\em resonance parameters}, and the indices $\lambda_1,\dots,\lambda_s$ where they first appear are the {\em resonances}. In the balance \eqref{balance0}, the resonance parameters $r_{n_{l-1}+1},\dots,r_{n_l}$ have the corresponding resonance $\lambda_l$. The number $m_l=n_l-n_{l-1}$ of resonance parameters of resonance $\lambda_l$ is the {\em multiplicity} of $\lambda_l$. We have
\[
n_l=m_0+m_1+\dots+m_l,
\]
where $m_0=1$ counts the resonance parameter $t_0$ of resonance $\lambda_0=-1$.

We certainly want the free parameters to really represent the variety of Laurent series solutions. This means that we should require that the {\em resonance vectors} of resonance $\lambda_l$ to form an $n\times m_l$ matrix of full rank $m_l$
\[
R_l=\begin{pmatrix}
\dfrac{\partial a_{1,\lambda_l}}{\partial r_{n_{l-1}+1}} & \cdots & \dfrac{\partial a_{1,\lambda_l}}{\partial r_{n_l}} \\
\vdots && \vdots \\
\dfrac{\partial a_{n,\lambda_l}}{\partial r_{n_{l-1}+1}} & \cdots & \dfrac{\partial a_{n,\lambda_l}}{\partial r_{n_l}}
\end{pmatrix}.
\]
For the resonance $\lambda_0=-1$, this means that the {\em basic resonance vector} $R_0=(-k_1a_{1,0},\dots,-k_na_{n,0})^T$ is non-zero. In fact, when we consider the independence among all the resonance parameters, we really should require that the {\em resonance matrix}
\[
R=(R_0\; R_1\; \cdots \; R_s)
\]
to have the full rank $n_s$.

The exact condition for Theorem A is the following.

\begin{definition*}
A balance \eqref{balance0} is {\em principal}, if $n_s=n$ and the resonance matrix is invertible.
\end{definition*}

\subsection{Indicial Normalization}

By substituting $t_0=t-(t-t_0)$, the coefficients $a_{i,j}$ become power series of $(t-t_0)$ with analytic functions of $t,r_2,\dots,r_n$ as coefficients. Further substituting these power series into \eqref{balance0}, the balance becomes
\begin{align}
u_i
&= a_{i,0}^{(0)}(t)(t-t_0)^{-k_i} +\dots +a_{i,\lambda_1-1}^{(0)}(t)(t-t_0)^{\lambda_1-1-k_i} \nonumber\\
&\quad +a_{i,\lambda_1}^{(0)}(t,r_2,\dots,r_{n_1})(t-t_0)^{\lambda_1-k_i} \nonumber\\
&\qquad +\dots +a_{i,\lambda_2-1}^{(0)}(t,r_2,\dots,r_{n_1})(t-t_0)^{\lambda_2-1-k_i}\nonumber\\
&\quad + \cdots  \label{balance1} \\
&\quad +a_{i,\lambda_s}^{(0)}(t,r_2,\dots,r_{n_s})(t-t_0)^{\lambda_s-k_i}+\cdots. \nonumber 
\end{align}
We have
\[
a_{i,\lambda_l}^{(0)}(t,r_2,\dots,r_{n_l})
=a_{i,\lambda_l}(t,r_2,\dots,r_{n_l})+c_{i,\lambda_l}(t,r_2,\dots,r_{n_{l-1}}).
\]
Therefore the new resonance matrix
\[
R^{(0)}=(R_0^{(0)}\; R_1^{(0)}\; \cdots \; R_s^{(0)}),\quad 
R_l^{(0)}=\begin{pmatrix}
\dfrac{\partial a_{1,\lambda_l}^{(0)}}{\partial r_{n_{l-1}+1}} & \cdots & \dfrac{\partial a_{1,\lambda_l}^{(0)}}{\partial r_{n_l}} \\
\vdots && \vdots \\
\dfrac{\partial a_{n,\lambda_l}^{(0)}}{\partial r_{n_{l-1}+1}} & \cdots & \dfrac{\partial a_{n,\lambda_l}^{(0)}}{\partial r_{n_l}}
\end{pmatrix},
\]
is obtained from the old one simply by replacing $t_0$ with $t$ 
\[
R^{(0)}(t,r_2,\dots,r_n)=R(t,r_2,\dots,r_n).
\]
The matrix is still invertible.

Up to rearranging the order of $u_i$, we may assume $k_1a_{1,0}^{(0)}(t)\ne 0$. Then we introduce a new variable $\tau$ by
\begin{equation}\label{change1}
u_1=\tau^{-k_1}.
\end{equation}
Fixing some $(-k_1)$-th root $\beta_1(t)$ of $a_{1,0}^{(0)}(t)$ and taking the $(-k_1)$-th root of the equality \eqref{balance1} for $u_1$, we get the power $\tau$-series
\begin{align}
\tau
&= \beta_1(t)(t-t_0) +\dots +\beta_{\lambda_1-1}(t)(t-t_0)^{\lambda_1-1} \nonumber \\
&\quad +\beta_{\lambda_1}(t,r_2,\dots,r_{n_1})(t-t_0)^{\lambda_1}+\dots+\beta_{\lambda_2-1}(t,r_2,\dots,r_{n_1})(t-t_0)^{\lambda_2-1}  \nonumber \\
&\quad + \cdots   \label{tauexpand} \\
&\quad +\beta_{\lambda_s}(t,r_2,\dots,r_{n_s})(t-t_0)^{\lambda_s}+\cdots, \nonumber 
\end{align}
We may reverse the power $\tau$-series and get
\begin{equation}\label{texpand}
t-t_0=b_1\tau+b_2\tau^2+b_3\tau^3+\cdots,
\end{equation}
where $b_i$ is an analytic function of the same variables as $\beta_i$. 

Substituting the power $\tau$-series \eqref{texpand} into the balance \eqref{balance1}, we get Laurent $\tau$-series
\begin{align}
u_i
&= a_{i,0}^{(0)}b_1^{-k_i}\tau^{-k_i}+[-k_ia_{i,0}^{(0)}b_1^{-1}b_2+a_{i,1}^{(0)}b_1]b_1^{-k_i}\tau^{1-k_i}+\cdots  \label{expand0}  \\
&\quad +[-k_ia_{i,0}^{(0)}b_1^{-1}b_{j+1}+a_{i,j}^{(0)}b_1^j+\text{terms involving }a_{i,<j}^{(0)}, b_{\le j}]b_1^{-k_i}\tau^{j-k_i}+\cdots.  \nonumber
\end{align}
For $u_1=\tau^{-k_1}$, the equality \eqref{expand0} gives the recursive relation that computes the coefficients $b_j$
\[
b_{j+1}=\dfrac{b_1^{j+1}}{k_1a_{1,0}^{(0)}}a_{1,j}^{(0)}+\text{terms involving }a_{1,<j}^{(0)}, b_{\le j}.
\]
Substituting the coefficients $b_j$ into the other $u_i$ in \eqref{expand0}, we get updated Laurent $\tau$-series
\begin{align}
u_i
&= a_{i,0}^{(1)}(t)\tau^{-k_i}+\dots +a_{i,\lambda_1-1}^{(1)}(t)\tau^{\lambda_1-1-k_i} \nonumber\\
&\quad +a_{i,\lambda_1}^{(1)}(t,r_2,\dots,r_{n_1})\tau^{\lambda_1-k_i} +\dots 
+a_{i,\lambda_2-1}^{(1)}(t,r_2,\dots,r_{n_1})\tau^{\lambda_2-1-k_i}\nonumber\\
&\quad + \cdots  \label{balance2} \\
&\quad +a_{i,\lambda_s}^{(1)}(t,r_2,\dots,r_{n_s})\tau^{\lambda_s-k_i}+\cdots, 
\qquad
i>n_0=1, \nonumber 
\end{align}
where
\[
a^{(1)}_{i,j}
=\left(a_{i,j}^{(0)}-\dfrac{k_ia_{i,0}^{(0)}}{k_1a_{1,0}^{(0)}} a_{1,j}^{(0)}\right)b_1^{j-k_i}
+\text{ terms involving }a_{1,<j}^{(0)}, a_{i,<j}^{(0)}, b_{\le j},\quad
i>1.
\]

For $j=\lambda_l$, $i>1$ and $n_{l-1}<s\le n_l$, the formula for $a^{(1)}_{i,j}$ implies
\[
\dfrac{\partial a^{(1)}_{i,\lambda_l}}{\partial r_p}
=\left(\dfrac{\partial a^{(0)}_{i,\lambda_l}}{\partial r_p}-\dfrac{k_ia_{i,0}^{(0)}}{k_1a_{1,0}^{(0)}} \dfrac{\partial a^{(0)}_{1,\lambda_l}}{\partial r_p}\right)b_1^{j-k_i}.
\]
Therefore the updated resonance matrix
\[
R^{(1)}=(R_1^{(1)}\; \cdots \; R_s^{(1)}),\quad 
R_l^{(1)}=\begin{pmatrix}
\dfrac{\partial a_{2,\lambda_l}^{(1)}}{\partial r_{n_{l-1}+1}} & \cdots & \dfrac{\partial a_{2,\lambda_l}^{(1)}}{\partial r_{n_l}} \\
\vdots && \vdots \\
\dfrac{\partial a_{n,\lambda_l}^{(1)}}{\partial r_{n_{l-1}+1}} & \cdots & \dfrac{\partial a_{n,\lambda_l}^{(1)}}{\partial r_{n_l}}
\end{pmatrix}
\]
is obtained by deleting the first column (which consists of entirely $0$) of the following matrix
\[
\begin{pmatrix}
b_1^{-k_2} & &  \\
& \ddots & \\
 & & b_1^{-k_n}
\end{pmatrix}
\begin{pmatrix}
-\dfrac{k_2a_{2,0}^{(0)}}{k_1a_{1,0}^{(0)}} & 1 && \\
\vdots & &\ddots & \\
-\dfrac{k_na_{n,0}^{(0)}}{k_1a_{1,0}^{(0)}} & && 1
\end{pmatrix}
R^{(0)}
\begin{pmatrix}
b_1^{\lambda_1}I_{m_1} & &  \\
& \ddots & \\
 & & b_1^{\lambda_s}I_{m_s}
\end{pmatrix}.
\]
In other words, up to multiplying the powers of $b_1$ to rows and columns, $R^{(1)}$ is obtained as part of the row operation that uses the first entry of $R^{(0)}$ to eliminate the other terms in the first column. In particular, the invertibility of $R^{(0)}$ implies the invertibility of $R^{(1)}$. The size is reduced by $1$ because the resonance parameter $t_0$ has been ``absorbed'' into the new variable $\tau$.

\subsection{Resonance Variable}

Next we introduce new variables to ``absorb'' the resonance parameters $r_{(1,n_1]}=(r_2,\dots,r_{n_1})$ of resonance $\lambda_1$. By rearranging $u_2,\dots,u_n$ if necessary, we may assume that the first $m_1\times m_1$ submatrix $A^{(1)}$ of $R^{(1)}$ is invertible. Then we have
\[
R^{(1)}
=\begin{pmatrix}A^{(1)} & C^{(1)} \\ B^{(1)} & D^{(1)}\end{pmatrix},
\]
where
\[
R_1^{(1)}
=\begin{pmatrix}A^{(1)} \\ B^{(1)} \end{pmatrix}, \quad
(R_2^{(1)}\; \cdots \; R_s^{(1)})
=\begin{pmatrix}C^{(1)} \\ D^{(1)} \end{pmatrix}.
\]

We introduce new variables $\rho_2,\dots,\rho_{n_1}$ by truncating the $\tau$-series of $u_2,\dots,u_{n_1}$ in \eqref{balance2} at $\tau^{\lambda_1-k_i}$
\begin{align}
u_i
&= a_{i,0}^{(1)}(t)\tau^{-k_i}+\dots +a_{i,\lambda_1-1}^{(1)}(t)\tau^{\lambda_1-1-k_i} \label{change2} \\
&\quad +a_{i,\lambda_1}^{(1)}(t,\rho_2,\dots,\rho_{n_1})\tau^{\lambda_1-k_i},
\quad 1<i\le n_1. \nonumber
\end{align}
Then we have the equalities
\begin{align*}
a_{i,\lambda_1}^{(1)}(t,\rho_2,\dots,\rho_{n_1})
&= a_{i,\lambda_1}^{(1)}(t,r_2,\dots,r_{n_1}) +\dots 
+a_{i,\lambda_2-1}^{(1)}(t,r_2,\dots,r_{n_1})\tau^{\lambda_2-1-\lambda_1} \\
&\quad +a_{i,\lambda_2}^{(1)}(t,r_2,\dots,r_{n_2})\tau^{\lambda_2-\lambda_1} +\cdots  \\
&\quad + \cdots  \\
&\quad +a_{i,\lambda_s}^{(1)}(t,r_2,\dots,r_{n_s})\tau^{\lambda_s-\lambda_1}+\cdots, 
\quad 1<i\le n_1.
\end{align*}
By the inverse function theorem, the invertibility of 
\[
A^{(1)}(t,r_2,\dots,r_{n_1})=\dfrac{\partial (a_{2,\lambda_1}^{(1)},\dots,a_{n_1,\lambda_1}^{(1)})}{\partial (r_2,\dots,r_{n_1})}
\]
means that the left side is a locally invertible map of $\rho_2,\dots,\rho_{n_1}$, and for small $\tau$, the right side is a locally invertible map of $r_2,\dots,r_{n_1}$. By taking the inverse of the left side map, we get the power $\tau$-series for the new variables
\begin{align*}
\rho_i
&= r_i +\beta_{i,1}(t,r_2,\dots,r_{n_1})\tau+\cdots \\
&\quad +\beta_{i,\lambda_2-\lambda_1}(t,r_2,\dots,r_{n_1},r_{n_1+1},\dots,r_{n_2})\tau^{\lambda_2-\lambda_1}+\cdots  \\
&\quad + \dots  \\
&\quad +\beta_{i,\lambda_s-\lambda_1}(t,r_2,\dots,r_{n_1},r_{n_1+1},\dots,r_{n_s})\tau^{\lambda_s-\lambda_1}+\cdots, 
\quad 1<i\le n_1.  
\end{align*}
By taking the inverse of the right side map, we get the power $\tau$-series for the resonance parameters
\begin{align}
r_i
&= \rho_i +b_{i,1}(t,\rho_2,\dots,\rho_{n_1})\tau+\cdots  \nonumber \\
&\quad +b_{i,\lambda_2-\lambda_1}(t,\rho_2,\dots,\rho_{n_1},r_{n_1+1},\dots,r_{n_1})\tau^{\lambda_2-\lambda_1}+\cdots  \nonumber \\
&\quad + \cdots   \label{rexpand1} \\
&\quad +b_{i,\lambda_s-\lambda_1}(t,\rho_2,\dots,\rho_{n_1},r_{n_1+1},\dots,r_{n_s})\tau^{\lambda_s-\lambda_1}+\cdots, 
\quad 1<i\le n_1.   \nonumber
\end{align}

Substituting \eqref{rexpand1} into the balance \eqref{balance2}, we get the Laurent $\tau$-series
\begin{align}
u_i
&= a_{i,0}^{(1)}(t)\tau^{-k_i}+\dots +a_{i,\lambda_1-1}^{(1)}(t)\tau^{\lambda_1-1-k_i}+a_{i,\lambda_1}^{(1)}(\rho_2,\dots,\rho_{n_1})\tau^{\lambda_1-k_i} \nonumber \\
&\quad +\left[\dfrac{\partial a_{i,\lambda_1}^{(1)}}{\partial (r_2,\dots,r_{n_1})}(t,\rho_2,\dots,\rho_{n_1})b_1+a_{i,\lambda_1+1}^{(1)}(t,\rho_2,\dots,\rho_{n_1})\right]\tau^{\lambda_1+1-k_i} \nonumber\\
&\quad +\cdots   \label{expand1} \\
&\quad +\Biggl[\dfrac{\partial a_{i,\lambda_1}^{(1)}}{\partial (r_2,\dots,r_{n_1})}(t,\rho_2,\dots,\rho_{n_1})b_j+a_{i,\lambda_1+j}^{(1)}(t,\rho_2,\dots,\rho_{n_1},r_{n_1+1},\dots) \nonumber \\
&\qquad\qquad +\text{ terms involving }a_{i,<\lambda_1+j}^{(1)},b_{<j}\Biggr]\tau^{\lambda_1+j-k_i}  +\cdots, \nonumber
\end{align}
where $i>1$ and $b_j=(b_{2,j},\dots,b_{n_1,j})^T$ is the vertical coefficient vector from \eqref{rexpand1}. For $1<i\le n_1$, comparing \eqref{expand1} with \eqref{change2} gives the recursive relation that computes the coefficients  $b_j$
\[
A^{(1)}(t,\rho_2,\dots,\rho_{n_1})b_j=-a_{(1,n_1],\lambda_1+j}^{(1)}(t,\rho_2,\dots,\rho_{n_1},r_{n_1+1},\dots)+\text{ lower terms},
\]
where 
\[
a_{(n',n],j}=(a_{n'+1,j},\dots,a_{n,j})^T.
\]
Then we substitute the formula for $b_j$ into the other $u_i$ in \eqref{expand1} and get their updated Laurent $\tau$-series
\begin{align}
u_i
&= a_{i,0}^{(2)}(t)\tau^{-k_i}+\dots +a_{i,\lambda_1-1}^{(2)}(t)\tau^{\lambda_1-1-k_i} \nonumber \\
&\quad +a_{i,\lambda_1}^{(2)}(t,\rho_2,\dots,\rho_{n_1})\tau^{\lambda_1-k_i}+\dots+a_{i,\lambda_2-1}^{(2)}(t,\rho_2,\dots,\rho_{n_1})\tau^{\lambda_2-1-k_i} \nonumber\\
&\quad +a_{i,\lambda_2}^{(2)}(t,\rho_2,\dots,\rho_{n_1},r_{n_1+1},\dots,r_{n_2})\tau^{\lambda_2-k_i} \nonumber\\
&\qquad +\dots 
+a_{i,\lambda_3-1}^{(2)}(t,\rho_2,\dots,\rho_{n_1},r_{n_1+1},\dots,r_{n_2})\tau^{\lambda_3-1-k_i} \nonumber\\
&\quad + \cdots  \label{balance3} \\
&\quad +a_{i,\lambda_s}^{(2)}(t,\rho_2,\dots,\rho_{n_1},r_{n_1+1},\dots,r_{n_s})\tau^{\lambda_s-k_i}+\cdots, 
\qquad
i>n_1, \nonumber 
\end{align}
where
\begin{align*}
a_{i,j}^{(2)}
&=a_{i,j}^{(1)},
\quad j\le \lambda_1, \\
a^{(2)}_{i,j}
&=a_{i,j}^{(1)}(t,\rho_2,\dots,\rho_{n_1},r_{n_1+1},\dots) \\
&\quad -\dfrac{\partial a_{i,\lambda_1}^{(1)}}{\partial (r_2,\dots,r_{n_1})}(t,\rho_2,\dots,\rho_{n_1})(A^{(1)})^{-1}a_{(1,n_1],j}^{(1)}(t,\rho_2,\dots,\rho_{n_1},r_{n_1+1},\dots) \\
&\quad + \text{ lower terms},\quad j>\lambda_1.
\end{align*}
Written in vector form, for $j>\lambda_1$ we have
\[
a_{(n_1,n],j}^{(2)}
=a_{(n_1,n],j}^{(1)}-B^{(1)}(A^{(1)})^{-1}a_{(1,n_1],j}^{(1)},\quad
B^{(1)}=\left.\dfrac{\partial a_{(n_1,n],\lambda_1}^{(1)}}{\partial (r_2,\dots,r_{n_1})}\right|_{r_i=\rho_i}.
\]

For $j=\lambda_l$, $l>1$, and $n_{l-1}<s\le n_l$, the formula for $a_{(n_1,n],j}^{(2)}$ implies that the updated resonance vectors of resonance $\lambda_l$ form the matrix
\begin{align*}
R_l^{(2)}
&=\dfrac{\partial a_{(n_1,n],\lambda_l}^{(2)}}{\partial (r_{n_{l-1}+1},\dots,r_{n_l})} \\
&=\dfrac{\partial a_{(n_1,n],\lambda_l}^{(1)}}{\partial (r_{n_{l-1}+1},\dots,r_{n_l})}
-B^{(1)}(A^{(1)})^{-1}\dfrac{\partial a^{(1)}_{(1,n_1],\lambda_l}}{\partial (r_{n_{l-1}+1},\dots,r_{n_l})} \\
&=(-B^{(1)}(A^{(1)})^{-1}\;\; I)R_l^{(1)}.
\end{align*}
Therefore the updated resonance matrix 
\[
R^{(2)}
=(R_2^{(2)}\; \cdots \; R_s^{(2)}) 
=(-B^{(1)}(A^{(1)})^{-1}\;\; I) \begin{pmatrix}C^{(1)} \\ D^{(1)} \end{pmatrix} 
=D^{(1)}-B^{(1)}(A^{(1)})^{-1}C^{(1)}
\]
is obtained as part of the row operation on $R^{(1)}$ that uses the invertible matrix $A^{(1)}$ to eliminate $B^{(1)}$
\[
\begin{pmatrix} I & O \\ -B^{(1)}(A^{(1)})^{-1} & I \end{pmatrix}R^{(1)}
=\begin{pmatrix} A^{(1)} & C^{(1)} \\ O & D^{(1)}-B^{(1)}(A^{(1)})^{-1}C^{(1)} \end{pmatrix}.
\]
In particular, the invertibility of $R^{(1)}$ implies the invertibility of $R^{(2)}$. 

The process continues and follows the same pattern. After introducing $\rho_{n_{l-2}+1},\dots,\rho_{n_{l-1}}$, we get an updated Laurent $\tau$-series for the remaining $u_i$, $i>n_{l-1}$. The coefficients of the series are functions of $t,\rho_2,\dots,\rho_{n_{l-1}}$ and $r_{n_{l-1}+1},\dots,r_n$. Moreover, we also have the resonance matrix $R^{(l-1)}$, which is an $(n-n_{l-1})\times (n-n_{l-1})$ invertible matrix. By rearranging the remaining $u_i$, we may assume that the first $m_l\times m_l$ matrix $A^{(l-1)}$ in $R^{(l-1)}$ is invertible. Then we introduce new variables by  truncating the Laurent $\tau$-series of $u_{n_{l-1}+1},\dots,u_{n_l}$ at the resonance $\lambda_l$ and replacing the resonance parameters with new variables $\rho_{n_{l-1}+1},\dots,\rho_{n_l}$.
\begin{align}
u_i
&= a_{i,0}^{(l)}(t)\tau^{-k_i}+\dots +a_{i,\lambda_1-1}^{(l)}(t)\tau^{\lambda_1-1-k_i} \nonumber \\
&\quad +a_{i,\lambda_1}^{(l)}(t,\rho_2,\dots,\rho_{n_1})\tau^{\lambda_1-k_i}+\cdots \nonumber \\
&\quad +\cdots \label{change3} \\
&\quad +a_{i,\lambda_{l-1}}^{(l)}(t,\rho_2,\dots,\rho_{n_{l-1}})\tau^{\lambda_{l-1}-k_i}+\cdots \nonumber \\
&\quad +a_{i,\lambda_l}^{(l)}(t,\rho_2,\dots,\rho_{n_{l-1}},\rho_{n_{l-1}+1},\dots,\rho_{n_l})\tau^{\lambda_l-k_i},
\quad n_{l-1}<i\le n_l. \nonumber
\end{align}
The invertibility of $A^{(l)}$ then implies that we can find the power $\tau$-series for the new variables
\begin{align}
\rho_i
&= r_i +\beta_{i,1}(t,\rho_2,\dots,\rho_{n_{l-1}},r_{n_{l-1}+1},\dots,r_{n_l})\tau+\cdots \nonumber \\
&\quad +\beta_{i,\lambda_{l+1}-\lambda_l}(t,\rho_2,\dots,\rho_{n_{l-1}},r_{n_{l-1}+1},\dots,r_{n_{l+1}})\tau^{\lambda_{l+1}-\lambda_l}+\cdots \nonumber \\
&\quad + \cdots   \label{rhoexpand2} \\
&\quad +\beta_{i,\lambda_s-\lambda_l}(t,\rho_2,\dots,\rho_{n_{l-1}},r_{n_{l-1}+1},\dots,r_{n_s})\tau^{\lambda_s-\lambda_l}+\cdots, 
\quad n_{l-1}<i\le n_l, \nonumber 
\end{align}
and the similar power $\tau$-series for the resonance parameters
\begin{align}
r_i
&= \rho_i +b_{i,1}(t,\rho_2,\dots,\rho_{n_l})\tau+\cdots  \nonumber \\
&\quad +b_{i,\lambda_{l+1}-\lambda_l}(t,\rho_2,\dots,\rho_{n_l},r_{n_l+1},\dots,r_{n_1})\tau^{\lambda_2-\lambda_1}+\cdots  \nonumber \\
&\quad + \cdots   \label{rexpand2} \\
&\quad +b_{i,\lambda_s-\lambda_l}(t,\rho_2,\dots,\rho_{n_l},r_{n_l+1},\dots,r_{n_s})\tau^{\lambda_s-\lambda_l}+\cdots, 
\quad n_{l-1}<i\le n_l.   \nonumber
\end{align}
Then we substitute \eqref{rexpand2} into the Laurent $\tau$-series of the remaining $u_i$, $i>n_l$. We get the updated Laurent $\tau$-series, in which the resonance parameters $r_{n_{l-1}+1},\dots,r_{n_l}$ are replaced by the new variables $\rho_{n_{l-1}+1},\dots,\rho_{n_l}$. The updated resonance matrix $R^{(l)}$ is obtained from the row operation on $R^{(l-1)}$ that eliminates the terms below $A^{(l-1)}$.

We get the whole change of variable after exhausting all resonances.

\subsection{Regularity}

Now we argue that the Laurent series in the balance \eqref{balance0} are converted to power series for the new variables $\tau,\rho_2,\dots,\rho_n$.

We have the power series \eqref{tauexpand} for $\tau$, with analytic functions of $t,r_2,\dots,r_n$ as coefficients. By taking the Taylor expansions in $(t-t_0)$ of the coefficients, we get the power series
\begin{align}
\tau
&= \alpha_1(t_0)(t-t_0) +\dots +\alpha_{\lambda_1-1}(t_0)(t-t_0)^{\lambda_1-1} \nonumber \\
&\quad +\alpha_{\lambda_1}(t_0,r_2,\dots,r_{n_1})(t-t_0)^{\lambda_1}+\cdots  \nonumber \\
&\quad + \cdots   \label{tauexpand3} \\
&\quad +\alpha_{\lambda_s}(t_0,r_2,\dots,r_{n_s})(t-t_0)^{\lambda_s}+\cdots, \nonumber 
\end{align}
with analytic functions of $t_0,r_2,\dots,r_n$ as coefficients. We note that 
$\alpha_1(t_0)=\beta_1(t_0)=a_{1,0}(t_0)^{-\frac{1}{k_1}}$ is nonzero.

We have the power series \eqref{rhoexpand2} for $\rho_i$. If $n_{l-1}<i\le n_l$, then the coefficients are analytic functions of $t$, $\rho_2,\dots,\rho_{n_{l-1}}$, $r_{n_{l-1}+1},\dots,r_n$. We may substitute the power series \eqref{tauexpand3} for $\tau$ and successively substitute the power series for $\rho_i$ corresponding to smaller resonances to the power series for $\rho_i$ corresponding to bigger resonances. Then we get updated power series for $\rho_i$, with analytic functions of $t,r_2,\dots,r_n$ as coefficients. Finally, we may take the Taylor expansions in $(t-t_0)$ of these coefficients and get the power series
\begin{align}
\rho_i
&= r_i +\alpha_{i,1}(t_0,r_2,\dots,r_{n_l})(t-t_0)+\cdots \nonumber \\
&\quad +\alpha_{i,\lambda_{l+1}-\lambda_l}(t_0,r_2,\dots,r_{n_{l+1}})(t-t_0)^{\lambda_{l+1}-\lambda_l}+\cdots \nonumber \\
&\quad + \cdots   \label{rhoexpand3} \\
&\quad +\alpha_{i,\lambda_s-\lambda_l}(t_0,r_2,\dots,r_{n_s})(t-t_0)^{\lambda_s-\lambda_l}+\cdots, 
\quad n_{l-1}<i\le n_l, \nonumber 
\end{align}

We conclude that the Laurent series for $u_1,\dots,u_n$ are converted to power series for $\tau,\rho_2,\dots,\rho_n$. Moreover, the new variables satisfy the initial conditions
\[
\tau(t_0)=0,\quad
\tau'(t_0)=a_{1,0}(t_0)^{-\frac{1}{k_1}}\ne 0,\quad
\rho(t_0)=r_i.
\]

Next we argue that the new system of differential equations for $\tau,\rho_2,\dots,\rho_n$ is regular. The formulae \eqref{change1} and \eqref{change3} for the change of variable imply that the new system is
\[
\tau'=g_1(t,\tau,\rho_2,\dots,\rho_n),\quad
\rho_i'=g_i(t,\tau,\rho_2,\dots,\rho_n),
\]
with the right side 
\begin{equation}\label{rightside}
g_i
=\gamma_i(t,\tau,\rho_2,\dots,\rho_n)
+\phi_{i,1}(t,\rho_2,\dots,\rho_n)\tau^{-1}
+\dots
+\phi_{i,N_i}(t,\rho_2,\dots,\rho_n)\tau^{-N_i},
\end{equation}
where $\gamma_i$ and $\phi_{i,j}$ are analytic in their variables. Since the Laurent series \eqref{balance0} is a formal solution of the original system, the power series \eqref{tauexpand3} and \eqref{rhoexpand3} is a formal solution of the new system. Substituting the power series into the $i$-th equation, we find that the left side is a power series, while the right side is a Laurent series, with $\alpha_1^{-N_i}\phi_{i,N_i}(t_0,r_2,\dots,r_n)(t-t_0)^{-N_i}$ as the lowest order term. Therefore we conclude that
\[
\phi_{i,N_i}(t_0,r_2,\dots,r_n)=0.
\]
Since this holds for all $t_0$ and resonance parameters $r_2,\dots,r_n$, and $\phi_{i,N_i}$ is analytic, we conclude that $\phi_{i,N_i}$ is constantly zero. This completes the proof of the regularity of the new system of equations.

\subsection{Triangular Change of Variable}

Our change of variable is not quite triangular. To get a triangular change of variable, we keep the first variable $\tau$ and modify the construction of new variables by introducing
\begin{align}
u_i
&= \tilde{a}_{i,0}^{(l)}(t)\tau^{-k_i}+\dots +\tilde{a}_{i,\lambda_1-1}^{(l)}(t)\tau^{\lambda_1-1-k_i} \nonumber  \\
&\quad +\tilde{a}_{i,\lambda_1}^{(l)}(t,\tilde{\rho}_2,\dots,\tilde{\rho}_{n_1})\tau^{\lambda_1-k_i}+\cdots  \nonumber \\
&\quad +\cdots \label{triangle} \\
&\quad +\tilde{a}_{i,\lambda_{l-1}}^{(l)}(t,\tilde{\rho}_2,\dots,\tilde{\rho}_{n_{l-1}})\tau^{\lambda_{l-1}-k_i}+\cdots  \nonumber \\
&\quad +\tilde{\rho}_i\tau^{\lambda_l-k_i},
\quad n_{l-1}<i\le n_l,  \nonumber
\end{align}
instead of \eqref{change3}.

The variables $\tilde{\rho}_i$ are related to $\rho_i$ by
\begin{align*}
\tilde{\rho}_i
&=a_{i,\lambda_l}^{(l)}(t,\rho_2,\dots,\rho_{n_{l-1}},\rho_{n_{l-1}+1},\dots,\rho_{n_l}) \\
&\quad +c_{i,\lambda_l}^{(l)}(t,\tau,\rho_2,\dots,\rho_{n_{l-1}}),
\quad n_{l-1}<i\le n_l,
\end{align*}
where $c_{i,\lambda_l}^{(l)}$ is analytic in its variables. The invertibility of $A^{(l)}$ implies that the relation is invertible, and the inverse is also analytic. The variables $\tilde{\rho}_i$ satisfy the initial condition
\[
\tilde{\rho}_{(n_{l-1},n_l]}(t_0)=A^{(l)} r_{(n_{l-1},n_l]}.
\]
Since $A^{(l)}$ are invertible, the initial values can be any number.

Finally, we remark that the proof for the regularity of the new system can be easily applied to general triangular change of variable
\begin{align*}
u_1 
&= \tau^{-k}, \\
u_i 
&= a_i(t,\tau,\rho_2,\dots,\rho_{i-1})+b_i(t,\tau,\rho_2,\dots,\rho_{i-1})\rho_i,\quad 1<i\le n.
\end{align*}
The key here is that, because $f_i$ are polynomial in $u_i$, the right side of the new system is of the form \eqref{rightside}. Then the power series solutions for $\rho_i$ with all the possible numbers as the initial values imply that $g_i$ have to be also analytic in $\tau$. 

The remark on the general triangular change of variable also shows that $f_i$ do not have to be polynomials. All we need is that the triangular change of variable converts the right side to be meromorphic in $\tau$. This often happens, for example, when $f_i$ are certain rational functions.

\section{Hamiltonian Structure} 
\label{hamilton}

Consider a Hamiltonian system 
\[
q'=\dfrac{\partial H}{\partial p},\quad
p'=-\dfrac{\partial H}{\partial q},
\]
where
\[
q=(q_1,\dots,q_n),\quad
p=(p_1,\dots,p_n),
\]
and $H=H(t,q,p)$ is analytic in $t$ and polynomial in $q$ and $p$. Consider a balance
\begin{align*}
q_i
&= a_{i,0}(t_0)(t-t_0)^{-k_i}+\dots
+a_{i,j}(t_0,r_2,\dots,r_{n_l})(t-t_0)^{j-l_i} +\cdots, \\
p_i
&= b_{i,0}(t_0)(t-t_0)^{-l_i}+\dots
+b_{i,j}(t_0,r_2,\dots,r_{n_l})(t-t_0)^{j-k_i} +\cdots,  
\end{align*}
where the coefficients for $(t-t_0)^{j-k_i}$ and $(t-t_0)^{j-l_i}$ depend only on the resonance parameters $r_2,\dots,r_{n_l}$ with resonance $\le j$.

The following is the exact condition for Theorem B. The definition is justified near the end of Section \ref{foundation}.

\begin{definition*}
A balance of the Hamiltonian system is {\em Hamiltonian principal}, if the resonance vectors form a simplectic basis of ${\mathbb R}^{2n}$, and there is $d$, such that $k_i+l_i=d-1$, and $\lambda_l+\mu_l=d-1$ for the resonances $\lambda_l$ and $\mu_l$ of simplectically conjugate resonance vectors.
\end{definition*}

The symplectic property depends on the order of vectors in the basis, and this order is not the same as the order of the column vectors for the resonance matrix. Consider the usual resonance matrix
\[
R=(R_0\; R_1\; \cdots \; R_{\bar{s}}),
\]
where the columns of $R_l$ have resonance $\lambda_l$. Since $\lambda_l$ is increasing in $l$, the condition $\lambda_l+\mu_l=d-1$ implies that a column vector in $R_l$ is simplectically conjugate to a column vector in $R_{\bar{s}-l}$. This is incompatible with our usual order in a simplectic basis $v_1,\dots,v_{2n}$, where $v_i$ is simplectically conjugate to $v_{n+i}$. So we need to reverse the last $n$ columns of $R$ and get the matrix
\[
S=R\begin{pmatrix}I_n & O \\ O & T_n \end{pmatrix},\quad
T_n=\begin{pmatrix}
0 & \cdots & 0 & 1 \\ 
0 & \cdots & 1 & 0 \\
\vdots && \vdots & \vdots \\
1 & \cdots & 0 & 0 
\end{pmatrix}.
\]
Then the definition requires that $S$ is a symplectic matrix
\[
S^TJS=J,\quad
J=\begin{pmatrix}O & I_n \\ -I_n & O \end{pmatrix}.
\]

Using the notation in the definition, we may denote the resonances by ($\bar{s}=2s+1$)
\[
-1=\lambda_0<\lambda_1<\dots<\lambda_s\le \mu_s<\dots<\mu_1<\mu_0=d,\quad
\mu_l=\lambda_{2s+1-l}.
\]
Since the resonance vectors form a simplectic basis, the columns of $R_l$ and $R_{2s+1-l}$ are simplectically conjugate and therefore the two matrices have the same number of columns. In case $\lambda_s=\mu_s$, the corresponding resonance vectors actually form the matrix $(R_s\;R_{s+1})$, where the columns of $R_s$ and $R_{s+1}$ are simplectically conjugate and therefore the two matrices have the same number of columns.

\subsection{Simplectic Resonance Variable}

The proof of Theorem A may be adapted to prove Theorem B. We only need to be more careful in keeping track of the simplectic structure.

Note that if the indicial normalization is assigned to $q_1$, then the last new variable, which corresponds to the largest resonance $\lambda_{2s+1}=\mu_0=d$, should be introduced for $p_1$. Therefore if we follow the construction in Section \ref{general}, we should rearrange the variables in the order $q_1,\dots,q_n$, $p_n,\dots,p_1$. The resonance matrix corresponding to this order is obtained by reversing the second half of the rows of $R$
\[
R'=\begin{pmatrix}I_n & O \\ O & T_n \end{pmatrix}R
=\begin{pmatrix}I_n & O \\ O & T_n \end{pmatrix}S\begin{pmatrix}I_n & O \\ O & T_n \end{pmatrix}.
\]

The resonance matrix is updated by row operations. This is equivalent to an $LU$ decomposition for $R'$, where $L$ is lower triangular and $U$ is upper triangular. In Section \ref{general}, the triangular shapes are actually blockwise. However, we get $LU$ decomposition only up to exchanging rows of the resonance matrix. For a Hamiltonian system, we need to exchange rows of $R'$ so that the simplectic structure is still preserved. 

We divide the simplectic matrix into four $n\times n$ blocks
\[
S=\begin{pmatrix} A & B \\ C & D \end{pmatrix},\quad
A^TC=C^TA,\quad
B^TD=D^TB,\quad
A^TD-C^TB=I_n.
\]
The columns of $\begin{pmatrix} A  \\ C \end{pmatrix}$ span a Lagrangian of the simplectic space ${\mathbb R}^{2n}$. It is well known in simplectic linear algebra that there is a subset $I\subset\{1,\dots,2n\}$, such that 
\begin{enumerate}
\item The projection $\pi_I\colon {\mathbb R}^{2n}\to {\mathbb R}^I$ sends the columns of $\begin{pmatrix} A  \\ C \end{pmatrix}$ to a basis of ${\mathbb R}^I$.
\item For any $1\le i\le n$, one of $i,n+i$ is in $I$, and one is not.
\end{enumerate}
Now for any $1\le i\le n$, if $n+i\in I$, then we exchange the $i$-th row and the $(n+i)$-th row of $S$, and if $n+i\in I$, then we do nothing. After this exchange, $S$ is still a simplectic matrix, with $A$ invertible. Correspondingly, this means that $q_i$ and $p_i$ are exchanged. But the system of differential equations is still Hamiltonian under such exchange.

So without loss of generality, we may assume that $A$ is invertible. By further  exchanging rows of $A$, we have $A=LU$. Exchanging rows of $A$ means exchanging $q_i$ with $q_j$. This should be balanced by exchanging $p_i$ with $p_j$ in order to keep the system Hamiltonian. Therefore if the $i$-th row and the $j$-th row of $A$ are exchanged, then the $i$-th row and the $j$-th row of $S$ should be exchanged, and the $(n+i)$-th row and the $(n+j)$-th row of $S$ should also be exchanged. Such exchange keeps $S$ to be simplectic.

After all the ``canonical exchanges'', which keep $S$ simplectic and the system Hamiltonian, we have a decomposition $A=LU$ that fits the construction in Section \ref{general}. Then the fact that $S$ is a simplectic matrix implies the $LU$ decomposition of the resonance matrix 
\[
R'=
\begin{pmatrix}
L & O \\ T_nCA^{-1}L & T_n(L^T)^{-1}T_n
\end{pmatrix}
\begin{pmatrix}
U & L^{-1}BT_n \\ O & T_n(U^T)^{-1}T_n
\end{pmatrix}.
\]
The construction of the change of variable in Section \ref{general} is guided by this $LU$ decomposition. The triangular version \eqref{change1} and \eqref{triangle} of the change of variable is
\begin{align}
q_1 &= \tilde{q}_1^{-k_1}, \nonumber \\
q_2 &= \alpha_{2,0}\tilde{q}_1^{-k_2}+\dots
+\alpha_{2,\lambda_1-1}\tilde{q}_1^{\lambda_1-1-k_2}
+\tilde{q}_2\tilde{q}_1^{\lambda_1-k_2}, \nonumber \\
&\vdots \nonumber \\
q_n &= \alpha_{n,0}\tilde{q}_1^{-k_n}+\dots
+\alpha_{n,\lambda_s-1}\tilde{q}_1^{\lambda_s-1-k_n}
+\tilde{q}_n\tilde{q}_1^{\lambda_s-k_n}, \nonumber \\
p_n &= \beta_{n,0}\tilde{q}_1^{-l_n}+\dots
+\beta_{n,\mu_s-1-l_n}\tilde{q}_1^{\mu_s-1-l_n}
+\tilde{p}_n\tilde{q}_1^{\mu_s-l_n}, \nonumber \\
&\vdots  \label{triangle2} \\
p_2 &= \beta_{2,0}\tilde{q}_1^{-l_2}+\dots
+\beta_{2,\mu_1-1-l_2}\tilde{q}_1^{\mu_1-1-l_2}
+\tilde{p}_2\tilde{q}_1^{\mu_1-l_2}, \nonumber \\
p_1 &= \beta_{1,0}\tilde{q}_1^{-l_1}+\dots
+\beta_{1,\mu_0-1-l_1}\tilde{q}_1^{\mu_1-1-l_1}
-k_1^{-1}\tilde{p}_1\tilde{q}_1^{\mu_0-l_1},  \nonumber 
\end{align}
where $\alpha_{i,j}$ is an analytic function of $t,\tilde{q}_2,\dots,\tilde{q}_{i-1}$, and $\beta_{i,j}$ is an analytic function of $t,\tilde{q}_2,\dots,\tilde{q}_n,\tilde{p}_n,\dots,\tilde{p}_{j+1}$. Note that a numerical factor $-k_1^{-1}$ is added in the last line.

\subsection{Canonical Change of Variable}

In this part, we argue that the change of variable is canonical
\[
dq_1\wedge dp_1+\dots+dq_n\wedge dp_n
=d\tilde{q}_1\wedge d\tilde{p}_1+\dots+d\tilde{q}_n\wedge d\tilde{p}_n.
\]

From the shape of the triangular change of variable \eqref{triangle2}, we have
\begin{align*}
dq_1 
&= -k_1\tilde{q}_1^{-k_1-1}d\tilde{q}_1
=-k_1\tilde{q}_1^{\lambda_0-k_1}d\tilde{q}_1, \\
dp_1 
&= -k_1^{-1}\tilde{q}_1^{\mu_0-l_1}d\tilde{q}_1
+\text{ linear combination of }\tilde{q}_1^{<\mu_0-l_1}, \\
dq_i
&= \tilde{q}_1^{\lambda_l-k_i}d\tilde{q}_i
+\text{ linear combination of }\tilde{q}_1^{<\lambda_l-k_i}, \\
dp_i 
&= \tilde{q}_1^{\mu_l-l_i}d\tilde{p}_i
+\text{ linear combination of }\tilde{q}_1^{<\mu_l-l_i}.
\end{align*}
By $(\lambda_l-k_i)+(\mu_l-l_i)=(\lambda_l+\mu_l)-(k_i+l_i) = (d-1)-(d-1)=0$, we have
\begin{align}
dq_1\wedge dp_1+\dots+dq_n\wedge dp_n
&=d\tilde{q}_1\wedge d\tilde{p}_1+\dots+d\tilde{q}_n\wedge d\tilde{p}_n \nonumber \\
&\quad +\tilde{q}_1^{-1}\omega_1+\dots+\tilde{q}_1^{-N}\omega_N, \label{form1}
\end{align}
where $\omega_i$ is a $2$-form in $\tilde{q}_1,\dots,\tilde{q}_n$, $\tilde{p}_1,\dots,\tilde{p}_n$, with analytic functions of $t,\tilde{q}_2,\dots,\tilde{q}_n$, $\tilde{p}_1,\dots,\tilde{p}_n$ as coefficients.

Think of the Laurent series in the Hamiltonian principal balance as a transform between $q_1,\dots,q_n$, $p_1,\dots,p_n$ and $t_0,r_2,\dots,r_{2n}$. We substitute the transform into the left side of \eqref{form1} and get a $2$-form in $t_0,r_2,\dots,r_{2n}$, with functions of $t,t_0,r_2,\dots,r_{2n}$ as coefficients. Since the symplectic form is invariant under the Hamiltonian flow, the coefficients are actually independent of $t$.

On the other hand, the Hamiltonian principal balance is converted to power series for the new variables $\tilde{q}_1,\dots,\tilde{q}_n$, $\tilde{p}_1,\dots,\tilde{p}_n$. Think of these power series as a transform between $\tilde{q}_1,\dots,\tilde{q}_n$, $\tilde{p}_1,\dots,\tilde{p}_n$ and $t_0,r_2,\dots,r_{2n}$. We substitute the transform into the right side of \eqref{form1} and get a $2$-form in $t_0,r_2,\dots,r_{2n}$, with functions of $t,t_0,r_2,\dots,r_{2n}$ as coefficients. Since this $2$-form is the same as the substitution of the Laurent series into the left side, the coefficients of this $2$-form are independent of $t$.

Let 
\[
\omega_N=d\tilde{q}_1\wedge\rho+\phi,
\]
where $\rho$ and $\phi$ are $1$-form and $2$-form in $\tilde{q}_2,\dots,\tilde{q}_n$, $\tilde{p}_1,\dots,\tilde{p}_n$, with analytic functions of $t$, $\tilde{q}_2,\dots,\tilde{q}_n$, $\tilde{p}_1,\dots,\tilde{p}_n$ as coefficients. After substituting the power series into $\omega_N$, we get
\[
\omega_N=dt_0\wedge \bar{\rho}+\bar{\phi}+(t-t_0)\bar{\psi},
\]
where
\begin{enumerate}
\item $\bar{\rho}$ and $\bar{\psi}$ are $1$-form and $2$-form in $r_2,\dots,r_{2n}$, with analytic functions of $t,t_0,r_2,\dots,r_{2n}$ as coefficients.
\item $\bar{\phi}$ is obtained as follows: Think of the initial data as a transform $(r_2,\dots,r_{2n})\mapsto (\tilde{q}_2(t_0),\dots,\tilde{q}_n(t_0), \tilde{p}_1(t_0),\dots,\tilde{p}_n(t_0))$ from the resonance parameters to new variables, with $t_0$ as a constant. Substituting the transform and $t=t_0$ into $\phi$ gives $\bar{\phi}$.
\end{enumerate}
The shape of $\omega_N$ implies that, after substituting the power series for $\tilde{q}_1,\dots,\tilde{q}_n$, $\tilde{p}_1,\dots,\tilde{p}_n$ into the right of \eqref{form1}, we get a $2$-form
\[
\alpha_1^{-N}(t-t_0)^{-N}\bar{\phi}+\cdots.
\]
Here $\alpha_1$ is the leading coefficient in \eqref{tauexpand3} for $\tau=\tilde{q}_1$, and depends only on $t_0$ and resonance parameters $r_i$ of resonance $0$. Moreover, $\cdots$ consists of terms that are either of order $>-N$ in $(t-t_0)$ or have $dt_0$ as a factor. Since we have already argued that the coefficients of the $2$-form should not depend on $t$, and $\alpha_1^{-N}$ is independent of $t$, and $\bar{\phi}$ is independent of $t$ and has no $dt_0$ factor, we conclude that $\bar{\phi}=0$. Since $\bar{\phi}$ is obtained from $\phi$ by an invertible transform, this implies $\phi=0$.

So we have
\[
\omega_N=d\tilde{q}_1\wedge\rho
=dt_0\wedge \bar{\rho}+(t-t_0)\bar{\psi}.
\]
Since there is no more $\phi$, we now know that $\bar{\rho}$ is obtained from $\rho$ in the similar way as $\bar{\phi}$ being obtained from $\phi$ above, and then multiplied by $-\alpha_1$. Moreover, $\bar{\psi}$ has the same description as above. Therefore the right side of \eqref{form1} now becomes
\[
\alpha_1^{-N}(t-t_0)^{-N}dt_0\wedge\bar{\rho}+\cdots,
\]
where $\cdots$ consists of terms that are of order $>-N$ in $(t-t_0)$. Since the coefficient of the $2$-form should not depend on $t$, and $\alpha_1$ and $\bar{\rho}$ are independent of $t$, we conclude that $\bar{\rho}=0$. Since $\bar{\rho}$ is obtained from $\rho$ by an invertible transform, this implies $\rho=0$.

We proved that $\omega_N=0$. Therefore there is no negative powers on the right of  \eqref{form1}, so that the change of variable \eqref{triangle2} is canonical.

The argument was carried out for the triangular change of variable. It also works for the blockwise triangular change of variable, making use of the simplectic property of the resonance matrix. Alternatively, we may verify (again by the simplectic property of the resonance matrix) that the transform between the new variables in the blockwise triangular change of variable and the new variables in the blockwise triangular change of variable is canonical.

\subsection{Change of the Hamiltonian Function}

In this part, we study the Hamiltonian function for the new Hamiltonian system.

Suppose the original Hamiltonian system is autonomous. Since our change of variable is canonical, the new system for the new variables is also an autonomous Hamiltonian system, and the new Hamiltonian function is simply obtained by applying the change of variable to the original Hamiltonian function. 

For an autonomous (not necessarily Hamiltonian) system, the coefficients in a balance do not depend on $t$. Then in the construction in Section \ref{general}, all the coefficients in the Laurent and power series are independent of $t$. Therefore the regularizing change of variable we get at the end is independent of $t$.

For the autonomous Hamiltonian system, applying the triangular change of variable \eqref{triangle2} to the original Hamiltonian function, we find the new Hamiltonian function to be of the form
\[
\tilde{H}(\tilde{q},\tilde{p})=h_0+h_1\tilde{q}_1^{-1}+\dots+h_N\tilde{q}_1^{-N},
\]
where $h_0$ is a polynomial of $\tilde{q}_1,\dots,\tilde{q}_n$, $\tilde{p}_1,\dots,\tilde{p}_n$, and $h_1,\dots,h_N$ are polynomials of $\tilde{q}_2,\dots,\tilde{q}_n$, $\tilde{p}_1,\dots,\tilde{p}_n$. Since the Hamiltonian function is a first integral of the Hamiltonian system, substituting the power series solution for $\tilde{q}_1,\dots,\tilde{q}_n$, $\tilde{p}_1,\dots,\tilde{p}_n$ into $\tilde{H}$ should give us an expression independent of $t$. The most singular term after the substitution is 
\[
h_N(\tilde{q}_2(t_0),\dots,\tilde{q}_n(t_0),\tilde{p}_1(t_0),\dots,\tilde{p}_n(t_0))\alpha_1^{-N}(t-t_0)^{-N}.
\]
Since $h_N$ and $\alpha_1$ are independent of $t$, we conclude that
\[
h_N(\tilde{q}_2(t_0),\dots,\tilde{q}_n(t_0),\tilde{p}_1(t_0),\dots,\tilde{p}_n(t_0))=0.
\]
Since the initial value map $(r_2,\dots,r_{2n})\mapsto (\tilde{q}_2(t_0),\dots,\tilde{q}_n(t_0),\tilde{p}_1(t_0),\dots,\tilde{p}_n(t_0))$ is invertible for triangular change of variable, we conclude that $h_N$ is a constantly zero function. This proves that $\tilde{H}(\tilde{q},\tilde{p})=h_0$ is a polynomial.

Now we turn to non-autonomous Hamiltonian system. In this case, the function $\tilde{H}(t,\tilde{q},\tilde{p})$ obtained by applying the time dependent change of variable to the original Hamiltonian function $H(t,q,p)$ is no longer gives the Hamiltonian function of the new system. In fact, $\tilde{H}$ may be singular in $\tilde{q}_1$. In what follows, we argue that it is the ``regular part'' $h_0(t,\tilde{q},\tilde{p})$ that becomes the Hamiltonian function of the new system.

We leave Hamiltonian systems for a moment and consider a general system of differential equations
\[
u'=f(t,u), \quad u=(u_1,\dots,u_n).
\]
For a change of variable $u=\varphi(t,\rho)$, we have
\[
u'=\dfrac{\partial \varphi}{\partial t}+\dfrac{\partial \varphi}{\partial \rho}\rho'.
\]
So the new system is
\[
\rho'=J_{\varphi}^{-1}f(t,\varphi(t,\rho))-J_{\varphi}^{-1}\dfrac{\partial \varphi}{\partial t},\quad
J_{\varphi}=\dfrac{\partial \varphi}{\partial \rho}.
\]

For the triangular change of variable \eqref{change1} and \eqref{triangle} ($\tilde{\rho}=(\tau,\tilde{\rho}_2,\dots,\tilde{\rho}_n)$ takes the place of $\rho$), we have 
\[
J_{\varphi}
=\begin{pmatrix}
-k_1\tau^{-1-k_1} & 0 & \cdots & 0 \\
\{\tau^{<\lambda_1-k_2}\} & \tau^{\lambda_1-k_2} & \cdots & 0 \\
\vdots & \vdots && \vdots \\
\{\tau^{<\lambda_s-k_n}\} & \{\tau^{<\lambda_s-k_n}\}& \cdots & \tau^{\lambda_s-k_n}
\end{pmatrix},
\]
where $\{\tau^{<j}\}$ means a (finite) linear combination of $\tau^{j'}$ with $j'<j$, with analytic functions of $t,\tilde{\rho}_2,\dots,\tilde{\rho}_n$ as coefficients. Moreover, $f(t,\varphi(t,\tilde{\rho}))$ and $\dfrac{\partial \varphi}{\partial t}$ are also (finite) linear combinations of $\tau^j$, with analytic functions of $t,\tilde{\rho}_2,\dots,\tilde{\rho}_n$ as coefficients. Then both $J_{\varphi}^{-1}f(t,\varphi(t,\tilde{\rho}))$ and $J_{\varphi}^{-1}\dfrac{\partial \varphi}{\partial t}$ are linear combinations of $\tau^j$ of the same type.

Let $J_{\varphi}^{-1}\dfrac{\partial \varphi}{\partial t}=(g_1,\dots,g_n)^T$. Then  $\dfrac{\partial \varphi}{\partial t}=J_{\varphi}(g_1,\dots,g_n)^T$, and by the formulae \eqref{change1} and \eqref{triangle}, we get
\begin{align*}
0 
& = -k_1\tau^{-1-k_1}g_1, \\
\{\tau^{<\lambda_1-k_2}\} 
& = \{\tau^{<\lambda_1-k_2}\}g_1+\tau^{\lambda_1-k_2}g_2 ,\\
& \vdots  \\
\{\tau^{<\lambda_s-k_n}\} 
& = \{\tau^{<\lambda_s-k_n}\}g_1+\{\tau^{<\lambda_s-k_n}\}g_2+\dots+\tau^{\lambda_s-k_n}g_n.
\end{align*}
By induction, it is easy to see that $g_1=0$, $g_2=\{\tau^{<0}\}$, $\dots$, $g_n=\{\tau^{<0}\}$. Therefore we conclude that $J_{\varphi}^{-1}\dfrac{\partial \varphi}{\partial t}$ is a linear combinations of $\tau^j$, $j<0$, with analytic functions of $t,\tilde{\rho}_2,\dots,\tilde{\rho}_n$ as coefficients. Since the differential equation 
\[
\tilde{\rho}'=J_{\varphi}^{-1}f(t,\varphi(t,\tilde{\rho}))-J_{\varphi}^{-1}\dfrac{\partial \varphi}{\partial t}
\]
obtained after the triangular change of variable \eqref{change1} and \eqref{triangle} is regular, we conclude that $J_{\varphi}^{-1}\dfrac{\partial \varphi}{\partial t}$ is the sum of terms in $J_{\varphi}^{-1}f(t,\varphi(t,\rho))$ with $\tau^j$, $j<0$. In other words, the new equation is 
\[
\tilde{\rho}'=[J_{\varphi}^{-1}f(t,\varphi(t,\rho))]_\text{regular},
\]
with the right side obtained by dropping terms with negative power of $\tau$.

Back to the non-autonomous Hamiltonian system, since the triangular change of variable \eqref{triangle2} is canonical, we have
\[
J_{\varphi}^{-1}\begin{pmatrix}\dfrac{\partial H}{\partial p} \\ -\dfrac{\partial H}{\partial q}\end{pmatrix}
=\begin{pmatrix}\dfrac{\partial \tilde{H}}{\partial \tilde{p}} \\ -\dfrac{\partial \tilde{H}}{\partial \tilde{q}}\end{pmatrix}.
\]
Here $\tilde{H}(t,\tilde{q},\tilde{p})$ is obtained by applying the triangular change of variable \eqref{triangle2} to $H(t,q,p)$ and is therefore a linear combination of $\tilde{q}_1^j$, with analytic functions of $t,\tilde{q}_2,\dots,\tilde{q}_n$, $\tilde{p}_1,\dots,\tilde{p}_n$ as coefficients. We separate the parts of $\tilde{H}$ with $\tilde{q}_1^j$, $j\ge 0$ and the rest part with $\tilde{q}_1^j$, $j<0$,
\[
\tilde{H}=[\tilde{H}]_\text{regular}+[\tilde{H}]_\text{singular}.
\]
Then it is easy to see that
\[
\dfrac{\partial \tilde{H}}{\partial \tilde{p}}
=\left[\dfrac{\partial \tilde{H}}{\partial \tilde{p}}\right]_\text{regular}
+\left[\dfrac{\partial \tilde{H}}{\partial \tilde{p}}\right]_\text{singular},
\]
with
\[
\left[\dfrac{\partial \tilde{H}}{\partial \tilde{p}}\right]_\text{regular}
=\dfrac{\partial [\tilde{H}]_\text{regular}}{\partial \tilde{p}},\quad
\left[\dfrac{\partial \tilde{H}}{\partial \tilde{p}}\right]_\text{singular}
=\dfrac{\partial [\tilde{H}]_\text{singular}}{\partial \tilde{p}},
\]
and the same happens to $\dfrac{\partial \tilde{H}}{\partial \tilde{q}}$. Thus we conclude that the new system of equations is
\[
\tilde{q}'=\dfrac{\partial [\tilde{H}]_\text{regular}}{\partial \tilde{p}},\quad
\tilde{p}'=-\dfrac{\partial [\tilde{H}]_\text{regular}}{\partial \tilde{q}}.
\]
This shows that $[\tilde{H}]_\text{regular}$ is the Hamiltonian function of the new system.

\section{A Mathematical Foundation for the Painlev\'e Test}
\label{foundation}

Let $f_i(t,u_1,\dots,u_n)$ be analytic in $t$ and polynomial in $u_1,\dots,u_n$. We consider a system of ordinary differential equations 
\begin{equation}\label{ode}
u_i'=f_i(t,u_1,\dots,u_n),\quad 1\le i\le n.
\end{equation}
The Painlev\'e test attempts to find all Laurent series solutions (or balances) with movable singularity $t_0$
\begin{equation}\label{balance}
u_i=c_i(t-t_0)^{-k_i}+a_{i,1}(t-t_0)^{1-k_i}+\dots+a_{i,j}(t-t_0)^{j-k_i}+\cdots.
\end{equation}
Here we use $c_i$ instead of $a_{i,0}$ to highlight the different role from the subsequent coefficients.

\subsection{Dominant Balance}

The first step in the Painlev\'e test is to determine the leading exponents $k_i$ and the leading coefficients $c_i$ of potential balances. There are usually several possible combinations of leading exponents and leading coefficients, which give several possible balances. 

\begin{definition*}
The leading exponents $k_1,\dots,k_n$ of a balance is {\em Fuchsian}, if the $k_*$-weighted degree of $f_i$ is $\leq k_i+1$.
\end{definition*}

The {\em $k_*$-weighted degree} of a polynomial in $u_1,\dots,u_n$ is obtained by taking the degree of $u_i$ to be $k_i$. 

If all the leading coefficients $c_i\ne 0$, then the choice of leading exponents is {\em natural}. For the natural leading exponents $k_*$, denote the dominant part of $f_i$
\[
D_i(t, u_1, \dots, u_n) 
=\sum \text{terms in $f_i$ with highest $k_*$-weighted degree}.
\]
The following is a simple criterion for the natural leading exponents to be Fuchsian. 

\begin{proposition}\label{prop1}
If $D_i(t_0,c_1,\dots,c_n)\ne 0$ and $k_ic_i\neq 0$ for each $i$, then the natural leading exponents is Fuchsian.
\end{proposition}

Let $d_i$ be the $k_*$-weighted degree of $f_i$. Then after substituting the balance \eqref{balance}, the lowest order term in $u_i'$ is $-k_ic_i(t-t_0)^{-k_i-1}$, and the lowest order term in $f_i(t,u_1,\dots, u_n)$ is
$D_i(t_0,c_1,\dots,c_n)(t-t_0)^{-d_i}$. Under the assumption of the proposition, both have nonzero coefficients. Therefore they are the ``real'' lowest order terms, and we conclude that
\[
d_i=k_i+1,\quad
D_i(t_0,c_1,\dots,c_n)=-k_ic_i.
\]
The first equality shows that the Fuchsian condition is satisfied.

If the natural exponents is not Fuchsian, then we need to choose somewhat unnatural leading exponents in order to satisfy the Fuchsian condition. For example, the Gelfand-Dikii hierarchy with 2 degrees of freedom is a Hamiltonian system with
\[
H = -q_1 p_2^2 - 2 p_1p_2  + 3 q_1^2q_2  - q_1^4 - q_2^2.
\]
One of the principal balances of the system is (we abbreviate $(t-t_0)$ as $t$ because the system is autonomous)
\begin{align*}
q_1 &
= t^{-2}+\frac{1}{3}r_2-\frac{1}{3}r_2^2t^2-\frac{2}{3}r_3t^3-\frac{10}{27}r_2^3t^4-\frac{1}{3}r_2r_3t^5-\frac{1}{3}r_4t^6+ \dots, \\
q_2 & 
=r_2t^{-2}-\frac{2}{3}r_2^2-r_3t-\frac{1}{3}r_2^3t^2+\left(-\frac{11}{54}r_2^4+\frac{3}{2}r_4\right)t^4+\cdots, \\
p_1 & 
=-t^{-5}+\frac{2}{3}r_2t^{-3}+\frac{1}{6}r_3-\frac{4}{27}r_2^3t-\frac{5}{6}r_2r_3t^2+\left(\frac{22}{81}r_2^4-\frac{11}{3}r_4\right)t^3+\cdots,  \\
p_2 & = 
t^{-3}+\frac{1}{3}r_2^2t+r_3t^2+\frac{20}{27}r_2^3t^3+\frac{5}{6}r_2r_3t^4+r_4t^5+ \dots.
\end{align*}
The natural leading exponents $2,2,5,3$ is not Fuchsian. By taking the leading exponent of 
\[
q_2=0t^{-4}+0t^{-3}+r_2t^{-2}-\frac{2}{3}r_2^2-r_3t-\frac{1}{3}r_2^3t^2+\left(-\frac{11}{54}r_2^4+\frac{3}{2}r_4\right)t^4+\cdots,
\]
to be $4$, we get leading exponents $2,4,5,3$ satisfying the Fuchsian condition. The corresponding leading coefficients are $1,0,-1,1$. The H\'enon-Heiles system in \cite{hy3} is another example.

When the Fuchsian condition is satisfied, we denote the dominant part of $f_i$
\[
f_i^D 
= \sum\text{terms in $f_i$ with $k_*$-weighted degree $k_i+1$}.
\]
Under the condition of Proposition \ref{prop1}, we have $f_i^D=D_i$. By substituting the balance \eqref{balance} into the system \eqref{ode} and comparing the coefficients of $(t-t_0)^{-k_i-1}$, we get
\begin{equation}\label{dominant}
f_i^D(t_0,c_1,\dots,c_n)=-k_ic_i.
\end{equation}
This is the (formal) {\em dominant balance} for the system \eqref{ode} at the (movable) singularity $t_0$.

\subsection{Kowalevskian Matrix}

Given the leading exponents and the leading coefficients, the next step of the Painlev\'e test is to find the subsequent coefficients in the balance. For $j>0$, by substituting the balance \eqref{balance} into the differential equations and comparing the coefficients of $(t-t_0)^{j-k_i-1}$, we get
\begin{align*}
(j-k_i)a_{i,j} 
&= \frac{\partial f_i^D}{\partial u_1}(t_0,c_1,\dots,c_n)a_{i,j}+\dots+\frac{\partial f_i^D}{\partial u_n}(t_0,c_1,\dots,c_n)a_{n,j} \\
& \quad +\text{ terms involving }c_*,a_{*,1},\dots,a_{*,j-1}.
\end{align*}
Define the {\em Kowalevskian matrix}
\begin{align*}\label{kmatrix}
K & = 
\begin{pmatrix}
\dfrac{\partial f_1^D}{\partial u_1}+k_1 & 
\dfrac{\partial f_1^D}{\partial u_2} & \cdots &
\dfrac{\partial f_1^D}{\partial u_n} \\
\dfrac{\partial f_2^D}{\partial u_1} & 
\dfrac{\partial f_2^D}{\partial u_2}+k_2 & \cdots &
\dfrac{\partial f_2^D}{\partial u_n} \\
\vdots & \vdots && \vdots \\
\dfrac{\partial f_n^D}{\partial u_1} & 
\dfrac{\partial f_n^D}{\partial u_2} & \cdots &
\dfrac{\partial f_n^D}{\partial u_n}+k_n
\end{pmatrix}  \\
& = 
\dfrac{\partial (f_1^D,\dots,f_n^D)}{\partial (u_1,\dots,u_n)}(t_0,c_1,\dots,c_n)
+\begin{pmatrix}
k_1 && \\ & \ddots & \\ && k_n
\end{pmatrix}.
\end{align*}
Then we have a recursive relation
\begin{equation}\label{recursion}
(K-jI)(a_{1,j},\dots,a_{n,j})^T = \text{ terms involving }c_*,a_{*,1},\dots,a_{*,j-1}.
\end{equation}

A {\em formal balance} is the Laurent series solution with leading coefficients obtained by solving \eqref{dominant} and with subsequent coefficients obtained by solving \eqref{recursion}.

We emphasize that the Fuchsian condition is needed for the definition of the Kowalevskian matrix. Ercolani and Siggia \cite{es} made a similar observation regarding the balance for the Gelfand-Dikii hierarchy. They suggested that one should not blindly proceed to use the linearized equation as the constraint on the resonance vectors and specifically pointed out that the leading exponent for $q_2$ should be changed from $2$ to $4$ in order to define the Kowalevskian matrix, which is
\[
K=\begin{pmatrix}
2 & 0 & 0 & -2 \\
-2 & 4 & -2 & -2 \\
12 & -6 & 5 & 2 \\
-6 & 2 & 0 & 3
\end{pmatrix}.
\] 

The solutions of \eqref{dominant} and \eqref{recursion} may not be unique. This leads to free parameters in the balance. Since \eqref{recursion} is a linear equation, we can choose the free parameters such that if a parameter does not appear in the leading coefficients, then the parameter first appears linearly. Although this property is not needed for our main theorems, it is a consequence of our set up.

There are three kinds of free parameters. The first kind is the free location $t_0$ of the singularity. The derivative of the balance \eqref{balance} in $t_0$ is 
\[
\frac{\partial u_i}{\partial t_0} 
= -k_ic_i(t-t_0)^{-k_i-1}+\left[(1-k_i)a_{1,i}+\dfrac{\partial c_i}{\partial t_0}\right](t-t_0)^{-k_i}+\cdots.
\]
Therefore the variation of the solution due to the parameter $t_0$ is characterized by the basic resonance vector $(-k_1c_1,\dots,-k_nc_n)$. 

Take the derivative of the system \eqref{ode} in $t$, we get
\[
u_i''
=\frac{\partial f_i}{\partial t}
+\frac{\partial f_i}{\partial u_1}u_1'+\dots
+\frac{\partial f_i}{\partial u_n}u_n'.
\]
Substituting in the balance \eqref{balance} and comparing the
coefficients of $(t-t_0)^{-k_i-2}$, we get
\[
(-k_i)(-k_i-1)c_i=\frac{\partial f_i^D}{\partial u_1}(t_0,c_1,\dots,c_n)(-k_1)c_1+\dots+\frac{\partial f_i^D}{\partial u_n}(t_0,c_1,\dots,c_n)(-k_n)c_n.
\]
The equalities can be rewritten as $(K+I)(-k_1c_1,\dots,-k_nc_n)=0$.

\begin{proposition}\label{resonance-1}
If the leading coefficients of a balance \eqref{balance} is Fuchsian, then the basic resonance vector $(-k_1c_1,\dots,-k_nc_n)$ is an eigenvector of the Kowalevskian matrix with eigenvalue $-1$.
\end{proposition}

The second kind of free parameters appear in the leading coefficients, and parameterize the subvariety given by \eqref{dominant}. The derivative of the equation \eqref{dominant} with respect to one such free parameter $r$ is
\[
-k_i\frac{\partial c_i}{\partial r}
=\frac{\partial f_i^D}{\partial u_1}\frac{\partial c_1}{\partial r}+\dots 
+\frac{\partial f_i^D}{\partial u_n}\frac{\partial c_n}{\partial r}.
\]
The equalities can be rewritten as $K\dfrac{\partial (c_1,\dots,c_n)^T}{\partial r}=0$.

\begin{proposition}\label{resonance0}
If the leading coefficients of a balance \eqref{balance} is Fuchsian, then the tangent vectors to the subvariety of allowable leading coefficients are the eigenvectors of the Kowalevskian matrix with eigenvalue $0$.
\end{proposition}

The third kind of free parameters appear in the subsequent coefficients. These are obtained from solving the recursive relation \eqref{recursion}. One has to worry about the compatibility between $K-jI$ and the right side of the recursive relation, which affects the existence of solutions. Leaving aside the existence issue for a moment, we may conclude the following from the recursive relation.

\begin{proposition}\label{resonancej}
If the leading coefficients of a balance \eqref{balance} is Fuchsian, then for $j>0$, the $j$-th coefficient vectors form an affine space parallel to the eigenspace of $K$ with eigenvalue $j$.
\end{proposition}

If $j$ is not an eigenvalue of $K$, then the $j$-th coefficients are uniquely determined by \eqref{recursion}. If $j$ is indeed an eigenvalue, then the $j$-th coefficients form a vector
\[
(a_{1,j},\dots,a_{n,j})^T=\hat{a}_j+r_1R_1+\dots+r_mR_m,
\]
where 
\begin{enumerate}
\item $\hat{a}_j$ depends only on $t_0,c_*,a_{*,1},\dots,a_{*,j-1}$.
\item $R_1,\dots,R_m$ form a basis of the eigenspace of $K$ with
eigenvalue $j$ and therefore depend only on $t_0$ and $c_*$.
\item $r_1,\dots,r_m$ can be any numbers.
\end{enumerate}

\begin{definition*}
A Fuchsian leading exponents is {\em principal}, if for all solutions of \eqref{dominant}, the following are satisfied.
\begin{enumerate}
\item The Kowalevskian matrix $K$ is diagonalizable, with the eigenspace of $-1$ being of dimension $1$, and all other eigenvalues being non-negative integers.
\item The recursive relation \eqref{recursion} is always consistent. 
\end{enumerate}
\end{definition*}

Given principal Fuchsian leading exponents, the corresponding principal balance \eqref{balance} is obtained by solving equations \eqref{dominant} and \eqref{recursion}. And the Painlev\'e test is passed.

The first condition implies that the eigenvalues of $K$ and their multiplicities are constants. By applying the implicit function theorem to \eqref{dominant}, we also know that the collection of leading coefficients form a submanifold with the eigenspace $\ker K$ of eigenvalue $0$ as the tangent space. Therefore the definition means exactly that the total number of free parameters, including $t_0$, is $n$. Moreover, the eigenvectors of $K$ form the resonance matrix.

For the balance of the Gelfand-Dikii hierarchy, the Kowalevskian matrix has eigenvalues $-1,2,5,8$, with the corresponding resonance matrix 
\[
R=\begin{pmatrix}
2 & 1 & -4 & -2 \\
0 & 3 & -6 & 9 \\
-5 & 2 & 1 & -22 \\
3 & 0 & 6 & 6
\end{pmatrix}.
\]

\subsection{Painlev\'e Test for Hamiltonian System}

Consider a Hamiltonian system given by a Hamiltonian function $H(t,q,p)$ that is analytic in $t$ and polynomial in $q$ and $p$. For a balance, we denote the leading exponents of $q_i$ and $p_i$ by $k_i$ and $l_i$.

\begin{definition}
A Hamiltonian system is {\em almost weighted homogeneous} relative to the leading exponents $k_1,\dots,k_n$, $l_1,\dots,l_n$ if $k_i+l_i=d-1$, where $d$ is the leading exponents weighted degree of $H$.
\end{definition}

The almost weighted homogeneous condition implies that the weighted degree of $\dfrac{\partial H}{\partial p_i}$ is $\le d-k_i=l_i+1$ and the similar inequality for $\dfrac{\partial H}{\partial q_i}$. Therefore the Fuchsian condition is satisfied. Let $H^D$ be the sum of those terms in $H$ with highest weighted degree and apply the same notation to other polynomial functions. Then $\left[\dfrac{\partial H}{\partial p_i}\right]^D=\dfrac{\partial H^D}{\partial p_i}$, and the Kowalevskian matrix is
\begin{align*}
K
&=	\begin{pmatrix}
	\dfrac{\partial^2 H^D}{\partial q_i\partial p_i}
	& \dfrac{\partial^2 H^D}{\partial p_i^2} \\
	-\dfrac{\partial^2 H^D}{\partial q_i^2} 
	& -\dfrac{\partial^2 H^D}{\partial p_i\partial q_i} 
  	\end{pmatrix}
+	\begin{pmatrix}
	k_1 &&& \\ & \ddots && \\ && l_1 & \\ &&& \ddots
	\end{pmatrix} \\
&=	\begin{pmatrix}
	O & I \\ -I & O  
	\end{pmatrix}
	\begin{pmatrix}
	\dfrac{\partial^2 H^D}{\partial q_i^2} 
	& \dfrac{\partial^2 H^D}{\partial p_i\partial q_i} \\
	\dfrac{\partial^2 H^D}{\partial q_i\partial p_i} 
	& \dfrac{\partial^2 H^D}{\partial p_i^2} 
	\end{pmatrix} 
+	\begin{pmatrix}
	k_1 &&& \\ & \ddots && \\ && l_1 & \\ &&& \ddots
	\end{pmatrix} \\
&=J{\mathcal H}+\Gamma.
\end{align*}
where ${\mathcal H}$ is the Hessian of $H^D$. The special shape of the
Kowalevskian matrix implies certain symplectic structures among the resonance
vectors.

\begin{proposition}
Suppose the Hamiltonian function is almost weighted homogeneous relative to the leading exponents of a balance. Suppose $v$ and $w$ are eigenvectors of the Kowalevskian matrix with eigenvalues $\lambda$ and $\mu$. If $\lambda+\mu\ne d-1$, then $\langle v,Jw \rangle=0$.
\end{proposition}

The condition $k_i+l_i=d-1$ means exactly $\Gamma J+J\Gamma=(d-1)J$. Then the lemma follows from
\begin{align*}
(\lambda +\mu )\langle v,Jw \rangle 
& = \langle Kv,Jw \rangle + \langle v,JKw \rangle \\
& = -\langle J(J{\cal H}+\Gamma)v,w \rangle
      +\langle v,J(J{\cal H}+\Gamma)w \rangle \\
& = \langle ({\cal H}-J\Gamma)v,w \rangle
      +\langle v,(-{\cal H}+J\Gamma)w \rangle \\
& = \langle v,\Gamma J w \rangle 
      +\langle v,J\Gamma w \rangle \\
& = (d-1)\langle v,Jw \rangle.
\end{align*}

The proposition implies that, by rescaling the vectors if necessary, we can find a simplectic basis of resonance vectors. For example, by rescaling and exchanging the last two columns of the resonance matrix for the balance of the Gelfand-Dikii hierarchy, we get a symplectic matrix
\[
R\begin{pmatrix}
1 & 0 & 0 & 0 \\
0 & 1 & 0 & 0 \\
0 & 0 & 0 & \frac{1}{9} \\
0 & 0 & -\frac{1}{81} & 0
\end{pmatrix}
=\begin{pmatrix}
2 & \frac{1}{3} & \frac{2}{81} & -\frac{4}{9} \\
0 & 1 & -\frac{1}{9} & -\frac{2}{3} \\
-5 & \frac{2}{3} & \frac{22}{81} & \frac{1}{9} \\
3 & 0 & -\frac{2}{27} & \frac{2}{3}
\end{pmatrix}.
\]
The last two columns are exchanged because the proposition tells us that in the original resonance matrix, up to rescaling, the first and fourth columns are simplectically dual to each other, and the second and third columns are simplectically dual to each other.

\end{document}